\documentclass[reqno]{amsart}
\usepackage{amssymb,url}
\usepackage{hyperref}
\date{\bf May 5, 2008}
\theoremstyle{plain}
\newtheorem{theorem}{Theorem}
\newtheorem{lemma}{Lemma}

\newtheorem{remark}{Remark}

\DeclareMathOperator{\dx}{d{\it x}\mspace{-2mu}}
\DeclareMathOperator{\dt }{d{\it t}\mspace{-2mu}}

\allowdisplaybreaks[4]

\begin{document}

\title[\hfilneg 05/2008\hfil ...New Proof on some sharp double integral Inequalities...]{New Proof on some sharp double integral Inequalities of the Hermite-Hadamard Type}

\author[V.N. Huy]{Vu Nhat Huy}
\address[V.N. Huy]{Department of Mathematics, Mechanics and Informatics\\
College of Science\\ Vi\d{\^{e}}t Nam National University\\ H\`{a} N\d{\^{o}}i, Vi\d{\^{e}}t Nam}

\author[W. J. Liu]{Wenjun Liu}
\address[W. J. Liu]{College of Mathematics and Physics\\
Nanjing University of Information Science and Technology \\
Nanjing 210044, China} \email{\href{mailto: W. J. Liu
<lwjboy@126.com>}{lwjboy@126.com}}

\author[Q. A. Ng\^{o}]{Qu\^{o}\hspace{-0.5ex}\llap{\raise 1ex\hbox{\'{}}}\hspace{0.5ex}c Anh Ng\^{o} $^\star$}
\address[Q. A. Ng\^{o}]{Department of Mathematics, Mechanics and Informatics\\
College of Science\\ Vi\d{\^{e}}t Nam National University\\ H\`{a} N\d{\^{o}}i, Vi\d{\^{e}}t Nam}
\email{\href{mailto: Q. A. Ng\^{o} <bookworm\_vn@yahoo.com>}{bookworm\_vn@yahoo.com}}



\begin{abstract}
In this paper, we derive a new proof on some sharp double integral inequalities of the Hermite-Hadamard type. Our approach is mainly based on well-known Taylor's theorem with the integral remainder.
\end{abstract}

\thanks{$^\star$ Corresponding Author}
\thanks{This paper was typeset using \AmS-\LaTeX}

\maketitle

\section{introduction}
Let $f(x)$ be a convex function on the closed interval $[a, b]$, the well-known Hermite-Hadamard's inequality  can be expressed as (\cite{dp}):
\begin{equation}\label{eq1}
f\left( {\frac{{a + b}}{2}} \right)
  \leqq  \frac{1}{b-a}\int\limits_a^b {f\left( x \right) \dx } \leqq
 \frac{{f\left( a \right) + f\left( b \right)}}{2}.
\end{equation}
It is well known that Hermite-Hadamard's inequality is an important cornerstone in mathematical analysis and optimization. There is a growing literature considering its refinements and interpolations now. Recently, Ujevi\'{c} obtained the following similar inequalities for convex functions
\begin{equation}\label{eq2}
 \frac{{f\left( a \right) + f\left( b \right)}}{2} -\frac{1}{8}S
  \leqq  \frac{1}{b-a}\int\limits_a^b {f\left( x \right) \dx } \leqq
  f\left( {\frac{{a + b}}{2}} \right)+\frac{1}{8}S,
\end{equation}
where $S=(f'(b)-f'(a))(b-a)$.

In this paper, we shall prove the following sharp double integral inequalities of the Hermite-Hadamard type.
\begin{theorem}\label{th1}
Let $I\subset \mathbb{R}$ be an open interval, $a, b\in I, a < b.$
If $f: I\rightarrow \mathbb{R}$ is   differentiable,   $m = \mathop
{\inf }\limits_{x\in\left[ {a,b} \right]} f'' \left( x \right)$ and
$M = \mathop {\sup }\limits_{x\in\left[ {a,b} \right]} f'' \left( x
\right)$. 
Then we have
\begin{equation}\label{eq3}
f\left( {\frac{{a + b}}{2}} \right) + \frac{m}{{24}}\left( {b - a}
\right)^2
  \leqq  \frac{1}{b-a}\int\limits_a^b {f\left( x \right) \dx } \leqq
 \frac{{f\left( a \right) + f\left( b \right)}}{2} - \frac{m}{{12}}\left( {b - a}
 \right)^2
\end{equation}
and
\begin{equation}\label{eq4}
 \frac{{f\left( a \right) + f\left( b \right)}}{2} - \frac{M}{{12}}\left( {b - a} \right)^2
  \leqq  \frac{1}{b-a}\int\limits_a^b {f\left( x \right) \dx } \leqq f\left( {\frac{{a + b}}{2}} \right) + \frac{M}{{24}}\left( {b - a}
\right)^2.
\end{equation}
The inequalities (\ref{eq3}) are sharp in the sense that the constants $\frac{1}{24}$ in the left-hand and $\frac{1}{12}$ in the right-hand cannot be replaced by a larger one, respectively. The inequalities (\ref{eq4}) are sharp in the sense that the constants $\frac{1}{12}$ in the left-hand and $\frac{1}{24}$ in the right-hand cannot be replaced by a smaller one, respectively.
\end{theorem}

\begin{remark}
{\rm (1)}\ If $f''\geqq 0,\, t\in [a, b]$, i.e. $f$ is a convex function, thus we can set $m = 0$ in (\ref{eq3}). Then, we recapture the well-known Hermite-Hadamard inequalities for convex functions.

{\rm (2)}\ If $f''\leqq 0,\, t\in [a, b]$, i.e. $f$ is a concave function, thus we can set $M = 0$ in (\ref{eq4}). Then, we get the following inequalities for concave functions.
\begin{equation}\label{eq5}
\frac{{f\left( a \right) + f\left( b \right)}}{2}   \leqq  \frac{1}{b-a}\int\limits_a^b {f\left( x \right) \dx } \leqq  f\left( {\frac{{a + b}}{2}} \right).
\end{equation}
\end{remark}

\begin{remark}
We also note that inequalities \eqref{eq3} and \eqref{eq4} have been proved in \cite{nu}. In this short note, we shall use an other approach.
\end{remark}

Before proving our main theorem, we need an essential lemma below. It is well-known in the literature as Taylor's formula or Taylor's theorem with the integral remainder.

\begin{lemma}[See \cite{AD}, Theorem 1]\label{bd}
Let $f :[a,b] \to \mathbb R$ and let $r$ be a positive integer. If $f$ is such that $f^{(r-1)}$ is absolutely continuous on $[a,b]$, $x_0 \in (a,b)$ then for all $x \in (a,b)$ we have
\[
f\left( x \right) = T_{r - 1} \left( {f,x_0 ,x} \right) + R_{r - 1} \left( {f,x_0 ,x} \right)
\]
where $T_{r - 1} \left( {f,x_0 , \cdot} \right)$ is Taylor's polynomial of degree $r-1$, that is,
\[
T_{r - 1} \left( {f,x_0 ,x} \right) = \sum\limits_{k = 0}^{r - 1} {\frac{{f^{\left( k \right)} \left( {x_0 } \right)\left( {x - x_0 } \right)^k }}
{{k!}}} 
\]
and the remainder can be given by
\begin{equation}\label{eq6}
R_{r - 1} \left( {f,x_0 ,x} \right) = \int\limits_{x_0 }^x {\frac{{\left( {x - t} \right)^{r - 1} f^{\left( r \right)} \left( t \right)}}
{{\left( {r - 1} \right)!}} \dt} 
\end{equation}
\end{lemma}

By a simple calculation, the remainder in \eqref{eq6} can be rewritten as
\[
R_{r - 1} \left( {f,x_0 ,x} \right) = \int\limits_0^{x - x_0 } {\frac{{\left( {x - x_0 - t} \right)^{r - 1} f^{\left( r \right)} \left( {x_0 + t} \right)}}{{\left( {r - 1} \right)!}} \dt} 
\]
which helps us to deduce a similar representation of $f$ as following
\begin{equation}\label{eq7}
f\left( {x + u} \right) = \sum\limits_{k = 0}^{r - 1} {\frac{{u^k }}{{k!}}f^{\left( k \right)} \left( x \right)} + \int\limits_0^u {\frac{{\left( {u - t} \right)^{r - 1} }}{{\left( {r - 1} \right)!}}f^{\left( r \right)} \left( {x + t} \right) \dt} .
\end{equation}

\section{Proofs of Theorem \ref{th1}}

Let
\[
F\left( x \right) = \int_a^x\limits {f\left( t \right) \dt } .
\]
Then
\[
F\left( b \right) = F\left( {\frac{{a + b}}
{2}} \right) + \frac{{b - a}}
{2}F'\left( {\frac{{a + b}}
{2}} \right) + \int\limits_0^{\frac{{b - a}}
{2}} {\left( {\frac{{b - a}}
{2} - t} \right)F''\left( {\frac{{a + b}}
{2} + t} \right) \dt } .
\]
Equivalently,
\[
F\left( b \right) = F\left( {\frac{{a + b}}
{2}} \right) + \frac{{b - a}}
{2}f\left( {\frac{{a + b}}
{2}} \right) + \int\limits_0^{\frac{{b - a}}
{2}} {\left( {\frac{{b - a}}
{2} - t} \right)f'\left( {\frac{{a + b}}
{2} + t} \right) \dt } .
\]
Similarly,
\begin{align*}
 F\left( a \right) = &F\left( {\frac{{a + b}}
{2}} \right) + \frac{{a - b}}
{2}f\left( {\frac{{a + b}}
{2}} \right) + \int\limits_0^{\frac{{a - b}}
{2}} {\left( {\frac{{a - b}}
{2} - t} \right)f'\left( {\frac{{a + b}}
{2} - t} \right) \dt } \hfill \\
 \mathop = \limits^{t: = - t} &F\left( {\frac{{a + b}}
{2}} \right) - \frac{{b - a}}
{2}f\left( {\frac{{a + b}}
{2}} \right) + \int\limits_0^{\frac{{b - a}}
{2}} {\left( {\frac{{b - a}}
{2} - t} \right)f'\left( {\frac{{a + b}}
{2} - t} \right) \dt }.
\end{align*}
Therefore,
\begin{align*}
 \int\limits_a^b {f\left( x \right) \dx } - \left( {b - a} \right)f\left( {\frac{{a + b}}
{2}} \right) &= F\left( b \right) - F\left( a \right) - \left( {b - a} \right)f\left( {\frac{{a + b}}
{2}} \right) \hfill \\
  &= \int\limits_0^{\frac{{b - a}}
{2}} {\left( {\frac{{b - a}}
{2} - t} \right)\left( {f'\left( {\frac{{a + b}}
{2} + t} \right) - f'\left( {\frac{{a + b}}
{2} - t} \right)} \right) \dt } \hfill \\
  &\geqq \int\limits_0^{\frac{{b - a}}
{2}} {\left( {\frac{{b - a}}
{2} - t} \right)2tm \dt } \hfill \\
  &= \frac{m}
{{24}}\left( {b - a} \right)^3 .
\end{align*}
On the other hand,
\begin{align*}
 \int\limits_a^b {f\left( x \right) \dx } & = F\left( b \right) - F\left( a \right) \hfill \\
  &= \left( {b - a} \right)F'\left( a \right) + \int\limits_0^{b - a} {\left( {b - a - t} \right)F''\left( {a + t} \right) \dt }
\end{align*}
and
\[
\frac{{b - a}}
{2}\left( {f\left( a \right) + f\left( b \right)} \right) = \frac{{b - a}}
{2}\left( {2f\left( a \right) + \int\limits_0^{b - a} {f'\left( {a + t} \right) \dt } } \right)
\]
which helps us to deduce that
\begin{align*}
 \frac{{b - a}}
{2}\left( {f\left( a \right) + f\left( b \right)} \right) &- \int\limits_a^b {f\left( x \right) \dx } \hfill \\
 & = \frac{{b - a}}
{2}\int\limits_0^{b - a} {f'\left( {a + t} \right) \dt } - \int\limits_0^{b - a} {\left( {b - a - t} \right)f'\left( {a + t} \right) \dt } \hfill \\
 & = \int\limits_0^{b - a} {\left( {t - \frac{{b - a}}
{2}} \right)f'\left( {a + t} \right) \dt } \hfill \\
 & = \int\limits_{\frac{{b - a}}
{2}}^{b - a} {\left( {t - \frac{{b - a}}
{2}} \right)f'\left( {a + t} \right) \dt } - \int\limits_0^{\frac{{b - a}}
{2}} {\left( {\frac{{b - a}}
{2} - t} \right)f'\left( {a + t} \right) \dt } \hfill \\
& = \int\limits_0^{\frac{{b - a}}
{2}} {\left( {\frac{{b - a}}
{2} - t} \right)f'\left( {b - t} \right) \dt } - \int\limits_0^{\frac{{b - a}}
{2}} {\left( {\frac{{b - a}}
{2} - t} \right)f'\left( {a + t} \right) \dt } \hfill \\
&  = \int\limits_0^{\frac{{b - a}}
{2}} {\left( {\frac{{b - a}}
{2} - t} \right)\left( {f'\left( {b - t} \right) - f'\left( {a + t} \right)} \right) \dt } \hfill \\
&  \geqq \int\limits_0^{\frac{{b - a}}
{2}} {\left( {\frac{{b - a}}
{2} - t} \right)\left( {b - a - 2t} \right)m \dt } \hfill \\
&  = \frac{m}
{{12}}\left( {b - a} \right)^3 .
\end{align*}

If we now substitute $f(x) = (x-a)^2$ in the inequalities
  then we find that
the left-hand side, middle term and right-hand side are all equal to
$\frac{(b-a)^2}{3}$. Thus, the inequalities (\ref{eq3}) are sharp in the
usual sense.

The proof of (\ref{eq3}) is completed.  The proof of (\ref{eq4}) is
similar.

\section*{Acknowledgements}

This work was supported by the Science Research Foundation of
Nanjing University of Information Science and Technology and the
Natural Science Foundation of Jiangsu Province Education Department
under Grant No.07KJD510133.

\end{document}